\documentclass{amsart}  
\newtheorem{theorem}{Theorem}
\newtheorem{lemma}{Lemma}
\newtheorem{proposition}{Proposition}
\newtheorem{corollary}{Corollary}
\newtheorem{example}{Example}
\newcommand{\C}{\mathbb{C}} 
\newcommand{\R}{\mathbb{R}}

\title[Approximation on arcs and dendrites]
{Approximation on arcs and dendrites going 
to infinity in $\C^n$ (extended version)} 
\thanks{Research supported by CRSNG(Canada), FCAR(Qu\'ebec) 
and Cinvestav(M\'exico)}
\date{\today}
\author[P. M. Gauthier E. S. Zeron]{P. M. Gauthier E. S. Zeron 
({\it In memoriam: Herbert James Alexander 1940-1999})}

\begin{document}

\begin{abstract}
The Stone-Weierstrass approximation theorem is extended to certain unbounded sets in
$\C^n.$ In particular, on a locally rectifiable arc going to infinity, each continuous
function can be approximated by entire functions.
\end{abstract}
\noindent

\maketitle\noindent
{\bf AMS subject classification numbers.} Primary: 32E30. Secondary: 32E25.\\
{\bf Key words:}Tangential approximation.  

\section{Introduction} 

This work is the original version of the paper: \textit{Approximation on
arcs and dendrites going to infinite in $\C^n$} \cite{GZ}. This version
could not be published in its extended form because of size limitations.
However, we wish to publish it because it contains a sketch of the proof
of Alexander-Stolzenberg's theorem, which we announced in \cite{GZ}, and
several lemas on tangential approximation by polynomial and meromorphic
functions which could not be included on \cite{GZ}. For example, we
include a not-very-know result of Arakelian in Proposition 5.

A famous theorem of Torsten Carleman \cite{C} asserts that for each continuous 
function $f$ on the real line $\R$ and for each positive continuous function 
$\epsilon $ on $\R,$ there exists an entire function $g$ on $\C$ such that 
$$|f(x)-g(x)|<\epsilon (x), \mbox{ for all } x\in \R.$$
Carleman's theorem was extended to $\C^n$ by Herbert Alexander \cite{Al1} who 
replaced the line $\R$ by a piecewise smooth arc going to infinity in $\C^n$ and by
Stephen Scheinberg \cite{Sc} who replaced the real line $\R$ by the real part $\R^n$
of $\C^n=\R^n+i\R^n.$
In the present work, we approximate on closed subsets of area zero in $\C^n$ 
and extend Alexander's theorem to locally rectifiable closed connected subsets
$\Gamma\subset\C^n$ which contain no closed curves.
 
Let $X$ be a subset of $\C^n.$ $X$ is a continuum if it is a compact 
connected set. The length and area of $X$ are the Hausdorff $1$-measure and 
$2$-measure of $X$ respectively. The set  $X$ is said to be of finite length at a point
$x\in X$ if this point has a neighbourhood in $X$ of finite length, and $X$ is said to
be of locally finite length if $X$ is of finite length at each of its points. Notice
that if $X$ is a set of locally finite length, then each compact subset of $X$ has
finite length (though $X$ itself need not be of finite length).  
We denote the polynomial hull of a compact set $X$ 
by $\widehat X$. The algebra of continuous functions defined on $X$ is denoted
by $\mathcal{C}(X)$. Finally, the definition and some properties of the 
first \v{C}ech cohomology group with integer coefficients $\check{H}^1(X)$ 
are presented in \cite{Gam} and \cite{St}.

\section{The Alexander-Stolzenberg theorem}

John Wermer laid the foundations of approximation on curves in $\C^n$ 
and prepared the way for a fundamental result of Gabriel  
Stolzenberg \cite{S1} concerning hulls and smooth curves 
(for history see \cite{S2}). In \cite{Al2}, Alexander comments that 
Stolzenberg's theorem can be improved 
to consider \textit{continua of finite length} instead of \textit{smooth curves}. 
We shall refer to the following version as the Alexander-Stolzenberg Theorem.

\begin{theorem}[Alexander-Stolzenberg] Let $X$ and $Y$ be two compact 
subsets of $\C^n,$ with $X$ polynomially convex and  $Y\setminus X$ 
of zero area. Then, 
\begin{description} 
\item[A] Every continuous function on $X\cup{Y}$ which is uniformly 
approximable on $X$ 
by polynomials is uniformly approximable on $X\cup{Y}$ by rational functions.
\end{description}
Suppose, moreover, there exists a continuum $\Upsilon\subset\C^n$ such that 
$\Upsilon\setminus{X}$ has locally finite length and $Y\subset(X\cup\Upsilon)$. 
Then:
\begin{description}
\item[B] $\widehat{X\cup{Y}}\setminus(X\cup{Y})$ is (if non-empty) a 
pure one-dimensional analytic subset of $\C^n\setminus(X\cup{Y}).$
\item[C] If the map 
$\check{H}^1(X\cup{Y})\rightarrow\check{H}^1(X)$ induced by 
$X\subset{X}\cup{Y}$ is injective, then $X\cup{Y}$ is polynomially convex.
\end{description}
\end{theorem}

The proof of this theorem is implicitly contained in the papers of 
Stolzenberg  \cite{S2} and Alexander \cite{Al2}, so we will devote this 
section to presenting a mere sketch of the proof by just indicating the 
necessary modifications to the existing proofs.  Note also, that in this 
Alexander-Stolzenberg Theorem, locally finite length is required only for 
parts \textbf{B} and \textbf{C}.   

The main arguments of the following lemma are essentially in 
\cite[p.~188]{S2}). 

\begin{lemma} Let $X$ and $Y$ be two compact subsets of $\C^n,$ with $X$ 
rationally convex and $Y\setminus{X}$ of zero area. Then, $X\cup Y$ is 
rationally convex. If, moreover, $X$ is polynomially convex, then 
given a point $p$ in the complement of $X\cup{Y}$, there is a
polynomial such that $f(p)=0$, $0\not\in{f}(X\cup{Y})$ and 
$\Re f(z)<-1$ for $z\in X$.
\end{lemma} 

\begin{proof} The set $X$ has a fundamental system of neighbourhoods which 
are rational polyhedra \cite[p. 283]{S1} or \cite{St}. Given a point $p$ in 
the complement of $X\cup{Y}$, choose a compact rational polyhedron 
$\widetilde{X}$ which contains $X$ in its interior, but 
$p\not\in\widetilde{X}$. Along with $\widetilde{X},$ the closure $K$ of 
$Y\setminus\widetilde{X}$ is also rationally convex because it has zero area 
\cite[p. 71]{Gam}, so there are two polynomials $g$ and $h$ such that 
$0\not\in{g}(K)$, $0\not\in{h}(\widetilde{X})$ and $g(p)=h(p)=0$. 

The rational function $(h/g)$ is smooth on $K$, and so $(h/g)(K)$ has zero 
area. Thus, we can find a complex number $\lambda\not\in{(h/g)}(K)$ whose 
absolute value $|\lambda|$ is so small that  the polynomial $f=h-\lambda{g}$ 
has no zeros on $\widetilde{X}\cup{K}$. Since 
$X\cup{Y}\subset\widetilde{X}\cup{K}$ and $f(p)=0$, 
it follows that $X\cup{Y}$ is rationally convex. 

If, in addition, $X$ is polynomially convex, one has just to choose 
$\widetilde{X}$ to be a compact polynomial polyhedron (see \cite{S1} or 
\cite[Lemma 7.4]{W}) and the polynomial $h$ to satisfy $\Re(h)<-1$ on 
$\widetilde{X}$; and so, for sufficiently small $\lambda$, $\Re(f)<-1$ on $X$. 
\end{proof}

To prove part \textbf{A} of Theorem 1, suppose first that the set $Y$ is 
itself of zero area. Stolzenberg's proof \cite[p.~187]{S2} uses the fact that 
the polynomial image of a finite union of smooth curves has zero area, and we 
also have that the image of $Y$ under a polynomial has area zero. Now suppose 
we merely know that $Y\setminus X$ has zero area. Let $f$ be a continuous 
function on $X\cup Y$ which is uniformly approximable on $X$ by polynomials 
and let $\epsilon > 0.$ There exists a polynomial $p$ such that 
$|f-p|<\epsilon /2$ on $X.$ Since $X$ is polynomially convex, it has a 
fundamental system of neighbourhoods which are polynomial polyhedra 
\cite[Lemma 7.4]{W}. From the continuity of $f-p,$ it follows that 
$|f-p|<\epsilon /2$ on some polynomial polyhedron $\widetilde X$ containing 
$X$ in its interior. Extend $p|_{\widetilde X}$ to a continuous function 
$\tilde p$ on $Y\setminus \widetilde X$ so that $|f-\tilde p|< \epsilon /2$ 
on $\widetilde X \cup Y.$ Since $\widetilde X \cup Y$ can be written as the 
union of $\widetilde X$ with a {\em compact} set of area zero, it follows 
from the first part of this proof that there is a rational function $h$ such 
that $|\tilde p-h| < \epsilon /2 $ on  $\widetilde X \cup Y.$ By the triangle 
inequality, $|f-h|<\epsilon $ on $X\cup Y$ which concludes the proof 
of \textbf{A}. 

Part \textbf{C} can be deduced from part \textbf{B} of Theorem 1 as presented 
in \cite[p.~188]{S2}. No changes are required because no metrical properties 
are invoked. Finally, the proof of part \textbf{B} is implicitly contained in 
Alexander paper \cite{Al2}, but we need to make several remarks.

Set $\Gamma=X\cup{Y}$ and suppose there is a point 
$p\in\widehat{\Gamma}\setminus\Gamma$. From Lemma 1, there is a polynomial 
$f$ such that $f(p)=0$, $0\not\in{f}(\Gamma)$ and $\Re(f)<-1$ on $X$. Fix the 
compact set $L=f(\Gamma)$ and the half-plane $H=\{(x,y)\in\C:x\geq-1/2\}$. 
Alexander's arguments \cite{Al2} can be slightly modified to show that 
$\widehat{\Gamma}\cap{f}^{-1}(\Omega)$ is a 1-dimensional analytic subset 
of $f^{-1}(\Omega)$, where $\Omega$ is the connected component of 
$\C\setminus L$ which contains the origin. Alexander uses the hypothesis that 
the set $L$ has finite length in the whole plane $\C$. However, his argument 
works even if we restrict the set $L$ to have finite length just in the 
half-plane $H$. Indeed, the intersection $L\cap{H}$ is the polynomial image 
of the compact set $\Gamma \cap{f}^{-1}(H)$ of finite length; recall that 
$\Gamma \cap{f}^{-1}(H)=(Y\setminus{X})\cap{f}^{-1}(H)$ has finite 
length because it is compact and contained in the  
set $\Upsilon\setminus{X}$ of locally finite length. Moreover, we shall see that we just 
need to rewrite Lemmas~3,~5 and~6 of \cite{Al2} to get the result.

For a set $X,$ let $\#X$ denote the number ($\leq\infty$) of elements in $X.$

\begin{lemma}[Lemma 3 of \cite{Al2}] Let $\Gamma$ be a compact set in $\C^n$ 
and  $f$ a polynomial in $\C^n$ such that $\Gamma\cap{f}^{-1}(H)$ has finite 
length. For $x\in\R$, set $N(x)=\#\{p\in\Gamma:\Re{f}(p)=x\}.$ Then 
$\int_{-1/2}^\infty{N}(x)dx<\infty$.
\end{lemma}  

\begin{lemma}[Lemma 5 of \cite{Al2}] Let $L\subset\C$ be compact 
and such that $\int_{-1/2}^\infty{N}(x)dx<\infty$ where 
$N(x)=\#\{q\in{L}:\Re(q)=x\}$. Then, for every component 
$\Omega$ of $\C\setminus L$ which meets the half-plane $H$, there exists a 
finite sequence $\Omega _0,\Omega _1,\ldots,\Omega _m$ of components of 
$\C\setminus L$ with $\Omega _0=$ the unbounded component, $\Omega _m=\Omega$ 
and $(\Omega _{j-1},\Omega _j)$ amply adjacent through rectangles  
$R_j=[a,b]\times[c_{j-1},c_j]$ contained in $H$ for $j=1,2,\ldots,m$.
\end{lemma}  

Lemma 2 need not be commented, and Lemma 3 holds by firstly choosing a 
horizontal line segment $[a,b]\times{c}\subset\Omega\cap{H}$ in the original 
proof.

\begin{lemma}[Lemma 6 of \cite{Al2}] Let $\Gamma$ be a compact subset of 
$\C^n$ and $f$ a polynomial in $\C^n$. Set $L=f(\Gamma)\subset\C$. Suppose 
that $\int_{-1/2}^\infty{N}(x)dx<\infty$ and that $L\cap{H}$ is contained in 
a continuum $L_1$ whose intersection $L_1\cap{H}$ is of finite length. Let 
$(\Omega_1,\Omega_2)$ be a pair of components of $\C\setminus(L\cup{L_1})$ 
which are amply adjacent through the square 
$R=[a,b]\times [c_1,c_2]\subset{H}$. Suppose 
$\widehat{\Gamma}\cap{f}^{-1}(\Omega_i)$ is a (possibly empty) pure 
1-dimensional analytic subset of $f^{-1}(\Omega_i)$ for $i=1$. Then, 
the same is true for $i=2$.
\end{lemma}

Alexander proves that $\widehat{\Gamma}\cap{f}^{-1}(D^o)$ is a pure 
1-dimensional analytic subset of $f^{-1}(D^o)$ where $D^o$ is an open set 
contained in $R\cap\Omega_2$. He deduces then that 
$\widehat{\Gamma}\cap{f}^{-1}(\Omega_2)$ is also a pure 1-dimensional
analytic set in $f^{-1}(\Omega_2)$ by using Lemma~11 of \cite{S2}. 
This lemma is quite amazing because the component $\Omega_2$ may 
not be completely contained in $H$. In our case, 
the analysis is done at a neighbourhood of the square $R\subset{H}$, so 
$\widehat{\Gamma}\cap{f}^{-1}(D^o)$ is analytic because $L_1$ has finite 
length in $H$. The result in Lemma 4 then follows, because 
the Stolzenberg Lemma 11 of \cite{S2}, which Alexander invokes, contains 
no metric restrictions in its hypotheses.

We conclude the proof of part \textbf{B} of Theorem~1 following Alexander's 
original arguments. If the equality $X\cup\Upsilon=\Gamma=X\cup{Y}$ holds, 
we let $\Omega$ be the connected component of $\C\setminus{L}$ which contains 
the origin. Apply Lemmas~2 and~3 to get a sequence 
$\Omega _0,\Omega _1,\ldots,\Omega _m=\Omega$. Finally, in Lemma 4, take 
$L_1=f(\Upsilon)\subset{L}$; recall that $\Upsilon\cap{f}^{-1}(H)$ is a 
compact set of finite length because it is contained in 
$\Upsilon\setminus{X}$. Then, noting that 
$f^{-1}(\Omega_0)\cap\widehat{\Gamma}=\emptyset$, we conclude inductively 
that $f^{-1}(\Omega)\cap\widehat{\Gamma}$ is a 1-dimensional analytic subset 
of $f^{-1}(\Omega)$. Hence: $\widehat{X\cup{Y}}\setminus(X\cup{Y})$ is 
analytic at an arbitrary point $p\in\widehat{X\cup{Y}}\setminus(X\cup{Y})$.

Now suppose that $X\cup{Y}$ is strictly contained in $X\cup\Upsilon$. Let 
$p\in\widehat{\Gamma}\setminus\Gamma$ as above. Modify $\Upsilon$ to obtain 
$\Upsilon_0$ such that $p\not\in\Upsilon_0$ but $\Upsilon_0$ is a continuum 
with $\Upsilon_0\setminus{X}$ of finite length and $Y\subset(X\cup\Upsilon_0)$
(say by radial projection to the boundary inside a ball containing $p$ in its 
interior, centered off $\Upsilon$, and disjoint from $\Gamma$). By the 
previous paragraph, $\widehat{X\cup\Upsilon_0}\setminus(X\cup\Upsilon_0)$ is 
analytic. By Lemma~7 of \cite{Al2}, $\widehat{X\cup{Y}}\setminus(X\cup{Y})$
is analytic at $p$. 

\vspace{10pt}
An arc $\Upsilon$, that is, the homeomorphic image of an interval of the real line, 
is of finite length at a point  $x\in \Upsilon $ if and only if $\Upsilon$ is 
(locally) rectifiable at $x.$
A direct consequence of the Alexander-Stolzenberg theorem is that every 
compact arc $\Upsilon$ which is locally rectifiable everywhere except perhaps
at finitely many of its points is polynomially convex and the approximation
condition $C(\Upsilon)=P(\Upsilon)$ holds; notice that $\Upsilon$ may be of 
infinite length. 

It is natural to ask whether the connectivity can be dropped in these 
considerations.  In fact, Alexander \cite{Al3} gave an example of a compact 
set $Y$ of finite length in $\C^2$ for which $\widehat Y\setminus Y$ is 
\textit{not} a pure one-dimensional analytic subset of $\C^2\setminus Y.$  
Thus, the connectivity cannot be dropped in the rectifiable Stolzenberg 
Theorem of Alexander.  Moreover, the following example shows that we cannot 
finesse Theorem 1 by enclosing $Y$ in a continuum of finite length,  although 
it is known that one can always construct a compact arc $\Gamma $ which 
meets every component of $Y$ (so $Y \cup \Gamma $ is connected) and 
$\Gamma \setminus Y$ is \textit{locally} rectifiable. 

\begin{example}
There exists a discrete bounded  set in $\C\setminus\{0\}$ such that 
no continuum containing this sequence has finite length. 
\end{example}

Consider the set $E$ consisting of the complex numbers 
$w_{j,k}=k/j^2+\sqrt{-1}/j$, for $j=1,2,\ldots$ and $k=0,1,\ldots,j$. It 
is easy to see that $E$ is contained in the disjoint union of the closed balls 
$\overline{B}_{j,k}$ with respective centers  $w_{j,k}$ and radii 
$\frac{1}{2(j+1)^2}$. Hence, each continuum which contains $E$ has to meet 
the center and the boundary of each ball $\overline{B}_{j,k}$, so its 
length has to be greater than $\sum_{j>1}\frac{j+1}{2(j+1)^2}=\infty.$

	\section{Approximation on unbounded sets}

We now pass from approximation on compacta to approximation on closed sets. 
Let $Y$ be a closed subset of $\C^n$ and $\mathcal{F}$ a subclass of 
$\mathcal{C}(Y).$  We say that a function $f$ defined on $Y$ can be 
uniformly (resp. tangentially) approximated by functions in $\mathcal{F}$ 
if for each positive constant $\epsilon $ (resp. positive continuous function
$\epsilon $ on $Y$) there is a $g\in \mathcal{F}$ such that $|f-g|<\epsilon$ 
on $Y.$ As $\mathcal{F}$ we are interested in the restrictions to $Y$ of the 
class $\mathcal{O}(\C^n)$ of entire functions and the class of meromorphic 
functions on $\C^n$ whose singularities do not meet $Y$. In the latter case,
we say that $f$ can be uniformly (resp. tangentially) approximated by 
meromorphic functions on $\C^n$. Recall that these meromorphic functions 
can be expressed as a quotient $p/q$ of entire functions $p$ and  $q$ with 
$q(z)\ne 0$ for all $z\in Y$ because the second Cousin problem can be solved 
in $\C^n$. 

If $Y$ is compact, then of course uniform and tangential approximation are 
equivalent and we may replace the classes of entire and meromorphic functions 
on $\C^n$ by the classes of polynomials and rational functions respectively. 
  
We say that $Y$ is a set of uniform (resp. tangential) approximation by 
functions in the class $\mathcal{F}$ if each $f\in \mathcal {C}(Y)$ can be 
uniformly (resp. tangentially) approximated  by functions in $\mathcal{F}.$
Of course, as we have defined them, such sets $Y$ cannot have any interior. 
In the literature, one also finds a more generous notion of sets of uniform or 
tangential approximation, which allows some sets having interior.  

Before going any further, we should point out that,  
sets of uniform approximation and sets of tangential approximation by holomorphic 
functions are in fact the same. This was proved by Norair Arakelian in  his doctoral 
dissertation \cite{Ar} in $\C.$ His proof works verbatim in $\C^n.$ 
Since this fact is not well known and the proof is short we include it. 

\begin{proposition}[Arakelian] Let $Y$ be a closed subset of $\C^n$ and let 
$\mathcal{F}$ be either the class of functions holomorphic on $Y$ or the 
class of entire functions. Then, $Y$ is a set of uniform approximation by 
functions in the class $\mathcal{F}$ if and only if it is a set of tangential 
approximation by functions in the same class. 
\end{proposition}

\begin{proof} Suppose $Y$ is a set of uniform approximation, 
$f\in\mathcal{C}(Y)$ and $\epsilon$ is a positive continuous function on 
$Y$. Set $\psi=\ln\epsilon.$ There exists a function $g_1\in\mathcal{F}$ such that 
$|\psi - g_1|<1$ on $Y.$  Setting $h = \exp (g_1-1),$ consider the functions 
$f/h \in C(Y).$ There exists a function $g_2\in\mathcal{F}$ such that $|f/h -g_2|<1$
on $Y.$  Then, $|f-hg_2| < |h|=\exp (\Re (g_1) -1) < \exp \psi = \epsilon .$ 
This completes the proof.  
\end{proof}

	The following is a non-compact version of the 
Stone-Weierstrass Theorem.

\begin{proposition}
A closed set $\Gamma\subset\C^n$ is a set of tangential approximation by 
entire functions if and only if one can approximate (in the tangential 
sense) the projections $\Re(z_m)$ for $m=1,\ldots,n$. 
\end{proposition}

\begin{proof} The necessity is trivial. Moreover, if one can 
approximate $\Re(z_m)$, one can approximate $\Im(z_m)$ as well 
since $\Im(z_m)=i(\Re(z_m)-z_m)$. Let $I$ be the natural diffeomorphism of 
$\C^n$ onto the real part $\R^{2n}$ of $\C^{2n}$. That is: $I_1(z)=\Re(z_1)$, 
$I_2(z)=\Im(z_1)$, $I_3(z)=\Re(z_2)$, $I_4(z)=\Im(z_2)$, etc., for $z\in\C^n$. 
Given two continuous function $f,\epsilon\in\mathcal{C}(\Gamma)$ with 
$\epsilon$ real positive, we may extend both of them continuously to all of 
$\C^n$ while keeping $\epsilon$ positive.  By the theorem of Scheinberg (see 
introduction), there is an entire function $F\in\mathcal{O}(\C^{2n})$ such that  
$|f(z)-F\circ{I}(z)|< \epsilon(z)/2$ for $z\in\C^n$.

Since $F$ is uniformly continuous on compact subsets of $\C^{2n},$ 
there is a positive continuous function $\delta $ on $\C^n$ such that 
$|F\circ{I}(z)-F(w)|<\epsilon (z)/2,$ for each $z\in \C^n$ and each
$w\in \C^{2n}$ for which $|I(z)-w|<\delta (z).$         

By hypotheses, we can approximate each $I_m$ on $\Gamma$ by entire 
functions and so there exists an entire mapping $h:\C^n\rightarrow \C^{2n}$ 
with $|I-h|<\delta$ on $\Gamma .$  Thus, $|F\circ I - F\circ h|<\epsilon/2$ 
on $\Gamma .$  By the triangle inequality, $|f-F\circ h| < \epsilon$ on 
$\Gamma .$ The function $F\circ h$ is entire because $h$ and $F$ are 
holomorphic.
\end{proof}

An interesting consequence of this result is that neither projection
$\Re $ nor $\Im $, in the complex plane $z\in\C$, can be
tangentially approximated  in the classical examples where the 
tangential approximation fails to hold, 
although uniform approximation may sometimes be possible.

\begin{example} Let
	$$Y = \bigcup_{j=0}^\infty Y_j,$$
where $Y_0 = [0,+\infty)\times \{0\},$ and for $j=1,2,\cdots ,$ 
$$Y_j = 
\left([0,j]\times \{\frac{1}{2j},\frac{1}{2j+1}\}\right) \cup \left(\{j\}\times [\frac{1}{2j},\frac{1}{2j+1}]\right).$$
Then, on $Y$ both $\Re $ and $\Im $ can be approximated uniformly but not tangentially by entire functions. 
\end{example}

\begin{proof} In his doctoral thesis, Arakelian \cite{Ar} gave a complete
characterization for sets of uniform approximation, from which it follows that $Y$ is
not a set of uniform approximation and {\it a fortiori} not a set of tangential
approximation.  Thus, by Proposition 2, $\Re $ and $\Im $ cannot be approximated
tangentially. We show that they can be approximated uniformly.

Fix $\epsilon > 0$ and set $Z_\epsilon = \{z:|\Im(z)|\le \epsilon \}$ and $W_\epsilon = Y\setminus  Z_\epsilon .$
We may assume $Z_\epsilon $ and  $W_\epsilon $ disjoint (by choosing an appropriate smaller $\epsilon $ if necessary). 
Now, define the function
$$f=\left\{\begin{array}{l}	\epsilon 	\mbox{ on } Z_\epsilon\\
				\Im 		\mbox{ on } W_\epsilon\end{array}\right..$$

Invoking again Arakelian's work (see \cite{Ar}, \cite[p.245]{GS} or \cite{Gai}), we deduce the existence of an entire
function $g$ such that $|f-g|<\epsilon $ on $Z_\epsilon \cup W_\epsilon .$  Hence, $|\Im -g|<2\epsilon $ on $Y,$ so $\Im $
and $\Re (z) = z-i\Im (z)$ can both be approximated uniformly on $Y$ by entire functions.   
\end{proof}

It is interesting to compare Propositions 1 and 2 in the light of the previous example. 

We should also notice that in Proposition 2 we can ask that the approximating
functions be holomorphic merely in a neighbourhood of $\Gamma$. 
We thus have that each continuous 
function $f\in\mathcal{C}(\Gamma)$ can be approximated (in the tangential 
sense) by functions holomorphic in a neighbourhood of $\Gamma$ if and only 
if every projection $\Re(z_m)$ can. This result suggests the following:

\begin{proposition}
Every closed set $\Gamma\subset\C^n$ of area zero is a set of tangential 
approximation by meromorphic functions in $\C^n.$ 
\end{proposition} 

\begin{proof} Let $f,\epsilon\in\mathcal{C}(\Gamma)$ be two continuous 
functions with $\epsilon$ real and positive, we must construct a 
meromorphic function $F$ such that $|F(z)-f(z)|<\epsilon(z)$ on $\Gamma$. 
Let $B_0$ be the empty set and $\overline{B}_k$ closed balls of radius 
$k$ and center in the origin.

\begin{lemma}
Each continuous function $h\in\mathcal{C}(\overline{B}_k\cup\Gamma)$ 
which can be uniformly approximated by polynomials in $\overline{B}_k$ 
can be uniformly approximated on 
$D=\overline{B}_k\cup(\Gamma\cap\overline{B}_{k+1})$ by rational 
functions whose singularities do not meet $\Gamma$.
\end{lemma}

\begin{proof}
From Theorem~1.A, there exists a rational function $a/b$ such that\\ 
$\left|(a/b)(z)-h(z)\right|<\delta$ for $z\in{D}$ and $0\not\in{b}(D)$. 
Notice that $b(\Gamma)$ has zero area, so we may choose a complex number 
$\lambda\not\in{b}(\Gamma)$ with absolute value so small such that 
$\lambda\not\in{b}(D)$ and 
$\left|\frac{a(z)}{b(z)-\lambda}-h(z)\right|<\delta$ for $z\in{D}$.
\end{proof}

The proof of the proposition now follows a classical inductive process. 
There exists  a rational function $F_1$ whose singularities do not meet 
$\Gamma$ and such that $|F_1(z)-f(z)|<(\frac{2}{3}-2^{-1})\epsilon(z)$ for 
$z\in\Gamma\cap\overline{B}_1$ by the previous lemma. Proceeding by 
induction, we shall construct a sequence of rational functions $F_k$ 
which converges uniformly on compact sets to a meromorphic function 
with the desired properties. 

Given a rational function $F_k$ whose singularities do not meet 
$\Gamma$ and such that $|F_k(z)-f(z)|<(\frac{2}{3}-2^{-k})\epsilon(z)$ 
in $\Gamma\cap\overline{B}_k$, let $h_k$ be a continuous function 
identically equal to zero on $\overline{B}_k$ and such that 
$|h_k(z)+F_k(z)-f(z)|<(\frac{2}{3}-2^{-k})\epsilon(z)$ for 
$z\in\Gamma\cap\overline{B}_{k+1}$ as well. Fix a real number $0<\lambda_k<1$ 
strictly less than $\epsilon(z)$ for every $z\in\Gamma\cap\overline{B}_{k+1}$. 

Applying Lemma~5, there exists a rational function $R_k$ whose 
singularities do not meet $\overline{B}_k\cup\Gamma$ and such that 
$|R_k(z)-h_k(z)|<2^{-1-k}\lambda_k$ for 
$z\in\overline{B}_k\cup(\Gamma\cap\overline{B}_{k+1})$. Thus, the 
singularities of the rational function $F_{k+1}(z)=F_k(z)+R_k(z)$ do not 
meet $\Gamma$ and $|F_{k+1}(z)-f(z)|<(\frac{2}{3}-2^{-1-k})\epsilon(z)$ 
for $z\in\Gamma\cap\overline{B}_{k+1}$ by the triangle  inequality.

Notice that $F_{k+1}(z)-F_k(z)$ is holomorphic and its absolute value 
is less than $2^{-1-k}$ inside $\overline{B}_k$, so the sequence $F_k$ 
converges to a meromorphic function with the desired properties.
\end{proof}

Similar inductive processes were originally employed to prove Carleman's 
theorem, stated in the introduction, which asserts that the real line $\R$ 
in $\C$ is a set of tangential approximation by entire functions. 
Alexander \cite{Al1} extended Carleman's theorem to piecewise smooth arcs 
$\Gamma$ going to infinity in $\C^n$. That is, $\Gamma$ is the the image 
of the real axis under a proper continuous embedding (a curve without 
self-intersections, \textit{going to infinity in both directions}). We should 
mention that this problem had been considered independently by Bernard 
Aupetit and Lee Stout (see Aupetit's book \cite{Au}). As a consequence of the 
Alexander-Stolzenberg Theorem, we also have the following further extension 
of Carleman's theorem, which was conjectured by Aupetit in \cite{Au} and 
announced by Alexander in \cite{Al1}. 

\begin{proposition}
Let $\Gamma $ be an arc which is locally rectifiable everywhere, except perhaps in a 
discrete subset, and going to infinity in $\C^n$. Besides, let $\epsilon$ be
a strictly positive continuous function on $\Gamma.$  Then, for each 
$f\in\mathcal{C}(\Gamma ),$ there exists an entire function $g$ on $\C^n$ 
such that $|f(z)-g(z)| < \epsilon (z), \mbox{ for all  } z\in \Gamma .$
That is, $\Gamma $ is a set of tangential approximation by entire functions. 
\end{proposition}

Alexander's proof (see also \cite{Au}), for the case that $\Gamma $ is smooth, 
relies ingeniously on the topology of arcs and the original Stolzenberg Theorem 
for smooth curves. 
It works also when the arc $\Gamma $ is locally rectifiable everywhere except
perhaps in a discrete subset. One only needs to rewrite Lemma~1 of \cite{Al1},
using the following corollary of Theorem~1.

\begin{corollary}
Let $X$ and $Y$ be two compact subsets of $\C^n$ such that $X$ is 
polynomially convex, $Y$ is connected and $Y\setminus{X}$ is locally 
of finite length everywhere except perhaps at finitely many of its points. 
If the map $\check{H}^1(X\cup{Y})\rightarrow\check{H}^1(X)$ induced by 
$X\subset{X}\cup{Y}$ is injective, then $X\cup{Y}$ is polynomially convex 
and every continuous function $f\in\mathcal{C}(X\cup{Y})$ which can be 
approximated by polynomials in $X$ can be approximated by polynomials 
on the union $X\cup{Y}$.
\end{corollary}

\begin{proof} Let $\{y_j\}$ be the points where $Y\setminus{X}$ is 
not of finite length. It is easy to see that $X\cup\{y_j\}$ is polynomially
convex and $f$ can be approximated by polynomials in $X\cup\{y_j\}$, so the 
result follows from Theorem~1, Lemma~1 and  the Oka-Weil theorem.
\end{proof}

We can also approximate by entire functions on unbounded sets  
which are more general than arcs, but first,  we need to
introduce the polynomially convex hull of non-compact sets:

\vspace{3mm} \textbf{Definition}. Given an arbitrary subset 
$Y$ of $\C^n$, its polynomially convex hull is defined by 
$\widehat{Y}=\bigcup\left\{\widehat{K}:K\subset{Y}\mbox{ is compact}\right\}$.

\begin{proposition} 
Let $\Gamma$ be a closed set in $\C^n$ of zero area such that
$\widehat{D\cup\Gamma}\setminus\Gamma$ is bounded for every compact set 
$D\subset\C^n$. Let $B_1$  be an open ball with center in the origin 
which contains the closure of $\widehat{\Gamma}\setminus\Gamma$. That is, 
the set $B_1\cup\Gamma$ contains the hull $\widehat{K}$ of every compact 
set $K\subset\Gamma$.

Then, given two continuous functions $f,\epsilon\in\mathcal{C}(\Gamma)$ 
such that $\epsilon$ is real positive and $f$ can be uniformly approximated 
by polynomials on $\Gamma\cap\overline{B}_1$, there exists an entire 
function $F$ such that $|F(z)-f(z)|<\epsilon(z)$ for $z\in\Gamma$.
\end{proposition} 

\begin{proof}
Let $B_0$ be the empty set, $B_1$ as in the hypotheses and $B_k$ open balls 
with center in the origin such that each $B_k$ contains the closure of 
$\widehat{\Gamma\cup\overline{B}}_{k-1}\setminus\Gamma$. That is, the set 
$B_k\cup\Gamma$ contains the hull $\widehat{K}$ of every compact set 
$K\subset(\Gamma\cup\overline{B}_{k-1})$. Define $X_k$ to be the polynomially 
convex hull of $\overline{B}_{k+1}\cap(\Gamma\cup\overline{B}_{k-1})$, so
$X_k\subset(B_k\cup\Gamma)$. The compact sets $X_k$ and  
$X_k\cap\overline{B}_k$ are both polynomially convex.

The given hypotheses automatically imply that there exists a polynomial
$F_1$ such that $|F_1(z)-f(z)|<(\frac{2}{3}-2^{-1})\epsilon(z)$ on 
$\Gamma\cap\overline{B}_1$. Proceeding by induction, we shall construct 
a sequence of polynomials $F_k$ which converges uniformly on compact 
sets to an entire function with the desired properties.

Given a polynomial $F_k$ such that 
$|F_k(z)-f(z)|<(\frac{2}{3}-2^{-k})\epsilon(z)$ on $\Gamma\cap\overline{B}_k$,
let $h_k$  be a continuous function equal to $F_k$ on $\overline{B}_k$ and 
such that $|h_k(z)-f(z)|<(\frac{2}{3}-2^{-k})\epsilon(z)$ for 
$z\in\Gamma\cap\overline{B}_{k+1}$ as well. Fix a real number $0<\lambda_k<1$ 
strictly less than $\epsilon(z)$ for every $z\in\Gamma\cap\overline{B}_{k+1}$.

Notice that $X_k=(X_k\cap\overline{B}_k)\cup(\Gamma\cap\overline{B}_{k+1})$. 
Hence, by Theorem 1.A, the function $h_k$ can be approximated by rational 
functions on $X_k$ because $X_k\cap\overline{B}_k$ is polynomially convex and 
$\Gamma$ has zero area. Moreover, the functions $h_k$  can be approximated by 
polynomials by the Oka-Weil theorem. Thus, there exists a polynomial $F_{k+1}$
such that $|F_{k+1}(z)-h_k(z)|<2^{-1-k}\lambda_k$ for $z\in{X}_k$, and so 
$|F_{k+1}(z)-f(z)|<(\frac{2}{3}-2^{-1-k})\epsilon(z)$ on 
$\Gamma\cap\overline{B}_{k+1}$.

Finally, the inequality $|F_{k+1}(z)-F_k(z)|<2^{-1-k}$ holds for 
$z\in\overline{B}_{k-1}$, so the sequence $F_k$ converges to an 
entire function with the desired properties.
\end{proof}

On the other hand, if the equality $\widehat{\Gamma}=\Gamma$ holds as well 
in the last proposition, we can choose the empty set instead of the open ball 
$B_1$  (because the proof is an inductive process); and so $\Gamma$ becomes 
a set of tangential approximation by entire functions. There are many closed 
sets $\Gamma$ which satisfy the hypotheses of the last proposition. For 
example, we have the following.

\begin{theorem}
Let $\Gamma$ be closed connected set of locally finite length in $\C^n$ whose 
first cohomology group $\check{H}^1(\Gamma)$ vanishes ($\Gamma$ contains no 
simple closed curves). Then, $\Gamma$ is a set of tangential approximation 
by entire functions.
\end{theorem}

\begin{proof} The proof strongly uses the topology of $\Gamma .$ We show that 
each point of $\Gamma $ has finite order, that is, has a basis of neighbourhoods in
$\Gamma$ having finite boundaries. Given a point $z\in\Gamma,$ let $B_r$ be the open 
ball in $\C^n$ 
of radius $r$ and center $z$. Since $\Gamma$ is locally of finite length, the 
intersection of $\Gamma$ with the closed ball $\overline B_r$ has finite 
length, so the intersection of $\Gamma$ with the boundary of $B_s$ must be a 
finite set for almost all radii $0<s<r.$ Whence, each sub-continuum of $\Gamma$ is 
locally connected \cite[p.~283]{Kr}. On the other hand, there are no simple closed 
curves contained in $\Gamma$ because $\check{H}^1(\Gamma)=0$, so each 
sub-continuum of $\Gamma$ is a dendrite, that is, a locally connected continuum
containing no simple closed curves.  
In particular, if $\Gamma$ is compact, then it is a dendrite. 

Notice the following lemma.

\begin{lemma}
Each compact subset $K\subset\Gamma$ is contained in a sub-continuum 
(dendrite) of $\Gamma$.
\end{lemma}

\begin{proof}
Since $\Gamma$ is locally connected, the set $K$ is contained in a 
finite union of sub-continua of $\Gamma$. The lemma now follows since 
$\Gamma$ is arcwise connected (see Theorem 3.17 of \cite{HY}).
\end{proof}

Let $D$ be a compact set in $\C^n$. Notice that $D\cup\Gamma$ may 
contain simple closed curves $\Upsilon$ with $D\cap\Upsilon\neq\emptyset$ 
but $\Upsilon\not\subset{D}$. We shall call such a simple closed curve
$\Upsilon\subset(D\cup\Gamma)$  a \textit{loop}. We show there exists a ball 
which contains all of these loops. Henceforth, let $B_r$ be open balls of 
radii $r$ and center in the origin, and choose a radius $s>0$ such that 
$D\subset{B}_s$. Recall that $\Gamma\cap\overline{B}_{s+1}$ has finite 
length, so there exists a ball $B_t$ with $s<t<s+1$ such that $\Gamma$ 
meets the boundary of $B_t$ only in a finite number of points 
$Q=\{q_1,\ldots,q_m\}$. Let $\{\Upsilon_j\}$ be the possible loops which 
meet the complement of $B_t$. The set $\bigcup\{\Upsilon_j\}\setminus{B}_t$ 
is contained in $\Gamma$ and can be expressed as the union of compact arcs 
(not necessarily disjoint) which lie outside of $\overline{B}_t$ except for
their two end points which lie in $Q$. Since $\Gamma $ cannot contain simple 
closed curves, two different arcs cannot share the same end points, and there 
can only be finitely many such arcs. Hence, there exists a ball $B_\delta$ 
which contains all the loops $\Upsilon,$ and $D\subset{B}_\delta$.

We shall show that $\widehat{D\cup\Gamma}\setminus\Gamma$ is bounded. 
Without loss of generality, we may suppose that $D$ is a closed ball. 
Since $\Gamma$ is connected, the hull $\widehat{D\cup\Gamma}$ is equal to 
$\bigcup_{r\geq\delta}\widehat{K}_r$, where $K_r$ is the connected component 
of $\overline{B}_r\cap(D\cup\Gamma)$ which contains $D$. We can prove that 
$\widehat{K}_r=\widehat{K}_\delta\cup{K}_r$, for every $r\geq\delta$, using 
Alexander's original argument. The following lemma is a literal translation 
of Lemma 1.(a) of \cite{Al1}, to our context.

\begin{lemma} For every $r\geq\delta$, 
$\widehat{K}_r=\widehat{K}_\delta\cup\tau_r$ 
where $\tau_r=\overline{K_r\setminus{K}_\delta}$.
\end{lemma}

Since the notation is quite complicated and different from Alexander's, and
we need to invoke Theorem 1.B, we shall include the proof of Lemma~7, but 
first we conclude the proof of the theorem. 

By Lemma 7, the set $\widehat{D\cup\Gamma}\setminus\Gamma$ is bounded because 
$\widehat{K}_r=\widehat{K}_\delta\cup\tau_r=\widehat{K}_\delta\cup{K}_r$ 
and $\widehat{D\cup\Gamma}=(\widehat{K}_\delta\cup\Gamma)$. Moreover, the
equality $\widehat{\Gamma}=\Gamma$ holds as well because each compact subset
of $\Gamma$ is contained in a dendrite of finite length and is polynomially 
convex (see Lemma~6 and Alexander's work \cite{Al2}), so we can deduce 
from Proposition~5 that $\Gamma$ is a set of tangential approximation.
\end{proof}

\begin{proof}[Proof of Lemma 7]
Let $T_r=\widehat{K}_\delta\cup\tau_r$ be the set on the 
right hand side of the asserted equality. Clearly, we have 
$T_r\subset\widehat{K}_r\subset\widehat{T}_r$ (the second inclusion is in 
fact equality). Thus it suffices to show that $T_r$ is polynomially convex. 
Arguing by contradiction, we suppose otherwise. By Theorem 1.B, 
$\widehat{T}_r\setminus{T}_r$ is a $1$-dimensional analytic subvariety of 
$\C^n\setminus{T}_r$. 

Let $V$ be a non-empty irreducible analytic component of 
$\widehat{T}_r\setminus{T}_r$. We claim that $\overline{V}\setminus{K}_r$ 
is an analytic subvariety of $\C^n\setminus{K}_r$. Since 
$T_r=\widehat{K}_\delta\cup\tau_r$, it suffices to verify 
this locally at a point $x\in\overline{V}\cap{Q}$ where
$$Q=\widehat{K}_\delta\setminus{K}_\delta.$$
By Theorem 1.B, both $\widehat{K}_r$ and $Q$ are analytic near $x$, where 
\textit{near} $x$ refers to the intersection of sets with \textit{small 
enough} neighbourhoods of $x$, here and below. Furthermore, near $x$, 
$\overline{V}\subset\widehat{K}_r$, $V\subset\widehat{K}_r\setminus{Q}$ and 
$Q\subset\widehat{K}_r$. Thus, near $x$, $Q$ is a union of some analytic 
components of $\widehat{K}_r$. It follows that near $x$, $\overline{V}$ is just a union 
of some of the other local analytic components of $\widehat{K}_r$ at $x$; in 
fact, near $x$, $\overline{V}=V\cup\{x\}$. Put
$$W=\overline{V}\setminus{K}_r.$$
Then $W$ is an irreducible analytic subset  of $\C^n\setminus{K}_r$ and moreover,
$$\overline{W}\setminus{W}\subset{K}_\delta\cup\tau_r=K_r.$$
Thus $\overline{W}\subset\widehat{K}_r$ by the maximum principle. 

Fix a point $p\in{V}\subset{W}$. Since $p\not\in{T}_r$, we have 
$p\not\in\widehat{K}_\delta$ and therefore there exists a polynomial $h$ such 
that $h(p)=0$ and $\Re{h}<0$ on $\widehat{K}_\delta$. By the open mapping 
theorem, either $h(W)$ is an open neighbourhood of $0$ or $h\equiv{0}$ on 
$W$. In the latter case, $h\equiv{0}$ on $\overline{W}$ and so 
$\overline{W}\setminus{W}$ is disjoint from $K_\delta$. This implies that
$\overline{W}\setminus{W}\subset\hat\tau_r$ so $W\subset\hat\tau_r$. We have
a contradiction because $\tau_r$ is contained in a  dendrite of finite length 
and is polynomially convex (see Lemma~6 and Alexander's work \cite{Al2}), and
moreover, a dendrite cannot contain a 1-dimensional analytic set. Hence, the 
former case holds. As $h(\tau_r)$ is nowhere dense in the plane (recall that 
it is of finite length), there is a small complex number $\alpha\in{h}(W)$ such 
that $\alpha\not\in{h}(\tau_r)$. Now put $g=h-\alpha$. If $\alpha$ is 
sufficiently small, we conclude that (i) $\Re{g}<0$ on $\widehat{K}_\delta$, 
(ii) $g(q)=0$ for some $q\in{W}$ and (iii) $0\not\in{g}(\tau_r)$.

Now (i) implies that the polynomial $g$ has a continuous logarithm on 
$\widehat{K}_\delta$ and so, by restriction, on $K_\delta$. We can extend this logarithm of
$g$ on $K_\delta$ to a continuous logarithm of $g$ on $K_r$ because of (iii), since the ball
$B_\delta $ was chosen such that every simple closed curve (loop) 
$\Upsilon\subset{K}_r$ is contained in $B_\delta $ and hence in $K_\delta .$
But $K_r$ contains $\overline{W}\setminus{W}$. Applying the argument 
principle \cite[p.~271]{S1} to $g$ on the analytic set $W$ gives a 
contradiction to (ii).
\end{proof}

\newpage

We remark that the condition of having zero area is
essential in Propositions 3 and 5, as the following example (inspired
by \cite{CGN}) shows.

\begin{example}
Let $\mathcal{I}$ be the closed unit interval $\left[0,1\right]$ 
of the real line and $K\subset\mathcal{I}$ the compact set 
$K=\left\{0,1,\frac{1}{2},\frac{1}{3},\frac{1}{4},\ldots\right\}$. 
It is easy to see that the $(2+\epsilon)$-dimensional Hausdorff measure 
of the closed connected set $Y=(\mathcal{I}\times\{0\})\cup(K\times\C)$ 
in $\C^2$ is equal to zero for every $\epsilon>0$, moreover, 
the equality $\widehat{Y}=Y$ holds. However, the following continuous 
function $f\in\mathcal{C}(Y)$ cannot be uniformly approximated by 
holomorphic functions in $\mathcal{O}(Y)$:
$$f(w,z)=\left\{\begin{array}{l}z\mbox{ if }
w=1\\0\mbox{ otherwise}\end{array}\right.$$
\end{example}

Suppose there exists a real number $\epsilon>0$ and a holomorphic function 
$g\in\mathcal{O}(Y)$ such that $|f-g|<\epsilon$ on $Y$. We automatically 
have that $g(w,z)$ is bounded, analytic and constant on each complex line 
$\{\frac{1}{j}\}\times\C$, $j=2,3,\dots$. Hence, the holomorphic function
$\frac{\partial{g}}{\partial{z}}$ vanishes on each complex line
$\{\frac{1}{j}\}\times\C$, $j=2,3,\dots$ as well. Since the zero set of 
$\frac{\partial{g}}{\partial{z}}$ is an analytic set, this derivative must 
be zero in a neighbourhood of $\{0\}\times\C$ and hence on the connected 
set $Y$. The last statement is a contradiction to the fact that 
$|g(1,z)-z|<\epsilon$ for every $z\in\C$.

On the other hand, to see that $\widehat{Y}=Y$, 
notice that $Y=\bigcup_{r>0}Y_r$, where 
$Y_r=(\mathcal{I}\times\{0\})\cup(K\times\Delta_r)$ and 
$\Delta_r\subset\C$ are closed discs of radius $r$. The set 
$K\times\Delta_r$ is polynomially convex because it is the Cartesian 
product of two polynomially convex sets in $\C$; and so $Y_r$ is
polynomially convex because of Theorem~1.

\vspace{10pt}
Although connectivity, as we have emphasized, plays a crucial role in this paper, similar
results can be obtained 
for sets whose connected components form a locally finite family. 
Finally, we remark that, on a Stein manifold, analogous results also hold by simply embedding
the Stein manifold into some $\C^n.$ A possible exception is Proposition 2, since $\Re (p)$
is not well-defined on a manifold. 

\newpage

\bibliographystyle{amsplain}

\noindent
{\bf Addresses}

\medskip
\noindent
D\'epartement de math\'ematiques et de statistique et \\
Centre de rech\`erches math\'ematiques, Universit\'e de Montr\'eal\\
Universit\'e de Montr\'eal, CP 6128 Centre Ville,\\
Montr\'eal, H3C 3J7, Canada\\
e-mail gauthier@ere.umontreal.ca

\smallskip
\noindent
Departamento de Matem\'aticas, Cinvestav I.P.N.\\
Apartado Postal 14-740, M\'exico D.F. 07000, M\'exico.\\
e-mail eszeron@math.cinvestav.mx

\end{document}